\documentclass[reqno]{amsart}
\usepackage{xcolor}

\usepackage{enumerate}
\usepackage{amsmath,amsthm,amscd,amsfonts,amssymb} 
\usepackage{latexsym}
\usepackage{mathrsfs}

\newtheorem{problem}{\bf Problem}
\newtheorem{cor}{Corollary}
\newtheorem{theorem}{Theorem}
\newtheorem{lemma}{Lemma}
\newtheorem*{conj}{\bf Conjecture}
\begin{document}

\title[A counter-example to CDSP]{An example of a cyclic analytic $2$-isometry with defect operator of rank $3$, whose Cauchy dual is not subnormal}


\author[S. A. Joshi]{Saee A. Joshi\textsuperscript{1}}
\address{\textsuperscript{1}Department of Mathematics, S.P. College, Pune, Maharashtra, India. 411 030\\ Email: saayalee@gmail.com}
\author[G. M. Phatak]{Geetanjali M. Phatak\textsuperscript{2}}
\address{\textsuperscript{2}Department of Mathematics, S.P. College, Pune, Maharashtra, India. 411 030\\ Email: gmphatak19@gmail.com}
\author[V. M. Sholapurkar]{Vinayak M. Sholapurkar\textsuperscript{3} }
\address{\textsuperscript{3}Bhaskaracharya Pratishthana, Pune, Maharashtra, India, 411004\\
	Email: vmshola@gmail.com}





\begin{abstract}
 
The Cauchy dual subnormality problem (CDSP, for short) asks whether the Cauchy dual of a $2-$isometry is subnormal. In this article, we provide a counter-example to CDSP by constructing a cyclic, analytic, $2-$isometry whose defect operator is of rank $3$. In particular, we prove that the Cauchy dual $M_z'$ of the multiplication operator $M_z$ on the Dirichlet space $D(\mu)$ is not subnormal if $\mu$ is supported at three equi-spaced points on the unit circle. 

\end{abstract}



\maketitle


\section{Introduction}\label{sec1}
The notion of Cauchy dual of an operator was introduced by Shimorin in \cite{Sh}.  The concept turns out to be interesting, especially while studying the interplay between two classes of  operators which are in some sense antithetical to each other. An excellent illustration of such classes is provided by contractive subnormal operators and completely hyperexpansive operators. Also, both the classes of operators have been described in terms of operator inequalities. The CDSP has been considerably dealt with in the literature. The problem has affirmative answer for some special classes of operators and negative answer in some other cases. For a detailed survey of these cases the reader is referred to \cite{cgr2022, mkvms2025, khasnis2025}.\\
We begin by quoting two examples, which led us to construct the counter-example presented in this article.
\begin{itemize}
	\item An example of a cyclic, analytic $2-$isometry $T$ such that $T^*T-I$ is of rank two and whose Cauchy dual is {\it subnormal} has been constructed in \cite[Example 7.2]{cgr2022} by choosing a measure on the unit circle supported at $\{1,-1\}$.
	\item An example of a cyclic, analytic $2-$isometry $T$ such that $T^*T-I$ is of rank two and whose Cauchy dual is {\it not subnormal} has been constructed in \cite{khasnis2025} by choosing a measure on the unit circle supported at $\{1,\zeta\}, \zeta\neq -1$.
\end{itemize}
The examples above suggest that the Cauchy dual $M_z'$ of $M_z$ on $D(\mu)$ would be subnormal if and only if measure $\mu$ is supported at the points which are equi-spaced on the unit circle.
However, in this article, we show that even if the measure $\mu$ is supported at $1,w,w^2$  where $w$ is the complex cube-root of unity, the Cauchy dual $M_z'$ of $M_z$ on $D(\mu)$ fails to be subnormal!\\The proof relies heavily on the techniques developed in \cite{cgr2022}. In particular, we capitalize on the fact that for a finitely supported measure $\mu$ on the unit circle the Dirichlet space $D(\mu)$ can be realized as a de Branges - Rovnyak space $\mathcal{H}(B)$ with the equality of norms \cite[Theorem 6.4]{cgr2022}. 

\section{Preliminaries}
Let $\mathbb{R}, \mathbb{C}, \mathbb{D}$ and $\mathbb{T}$ respectively denote the set of real numbers, the set of complex numbers, the open unit disc and the unit circle.

The notion of a reproducing kernel Hilbert space (RKHS) is well known and has been extensively studied, explored and crucially used in the theory of operators on Hilbert spaces. We quickly recall the definition of RKHS for a ready reference. \\
Let $S$ be a non-empty set and $\mathcal{F}$$(S,\mathbb{C})$ denote set of all functions from $S$ to $\mathbb{C}$. We say that $\mathcal{H}\subseteq \mathcal{F}$$(S,\mathbb{C})$ is a {\it Reproducing Kernel Hilbert Space}, briefly {\it RKHS}, if
\begin{enumerate}
\item[i.] $\mathcal{H}$ is a vector subspace of $\mathcal{F}$$(S,\mathbb{C})$
\item[ii.] $\mathcal{H}$ is equipped with an inner product $\langle , \rangle$ with respect to which $\mathcal{H}$ is a Hilbert space
\item[iii.] For each $x\in S$, the linear evaluation functional $E_x:\mathcal{H}$$\to \mathbb{C}$ defined by $E_x(f)=f(x)$, is bounded.
\end{enumerate}
If $\mathcal{H}$ is an RKHS on $S$ then for each $x\in S$, by Riesz representation theorem, there exists unique vector $k_x\in \mathcal{H}$ such that, for each $f\in \mathcal{H}$, $f(x)=E_x(f)=\langle f,k_x\rangle$.\\
The function $K:S\times S \to \mathbb{C}$ defined by $K(x,y)=k_y(x)$ is called as the reproducing kernel for $\mathcal{H}$. A classical reference for reproducing kernel Hilbert spaces is \cite{paulsen2016RKHS}.\\

We now recall the definitions of three special types of RKHS which are relevant in the present context. 
\begin{enumerate}

\item {\bf Hardy Space:} Let Hol$(\mathbb{D})$ denote the set of all holomorphic functions on $\mathbb{D}$. The Hardy space $H^2$ is defined as
\[H^2=\Bigg\{f(z)=\displaystyle\sum_{n=0}^{\infty}a_nz^n\in \text{Hol}(\mathbb{D}) : \displaystyle\sum_{n=0}^{\infty}|a_n|^2<\infty\Bigg\}.\]
It is a well known fact that $H^2$ is a reproducing kernel Hilbert space with the inner product defined as $\langle f,g \rangle :=\displaystyle\sum_{n=0}^{\infty}a_n\overline{b_n}$, if $f(z)=\displaystyle\sum_{n=0}^{\infty}a_nz^n$,\\$g(z)=\displaystyle\sum_{n=0}^{\infty}b_nz^n$ and a kernel function(namely, {\it Szego Kernel})\\ $K_\lambda(z)=\displaystyle\frac{1}{1-\overline{\lambda}z}$. The interested reader may consult \cite{douglas2012}.\\

\item {\bf Dirichlet spaces : }  For a finite positive Borel measure $\mu$ on unit circle $\mathbb{T}$, the Dirichlet space $D(\mu)$ is defined by
\[D(\mu):=\Bigg\{f\in \text{Hol}(\mathbb{D}) : \displaystyle\int_{\mathbb{D}}|f'(z)|^2P_\mu(z)dA(z)<\infty \Bigg\}.\]
Here, $P_\mu(z)$ is the Poisson integral for measure $\mu$ given by  $\displaystyle\int_{\mathbb{T}}\displaystyle\frac{1-|z|^2}{|z-\zeta|^2}d\mu(\zeta)$,\\ \noindent$dA$ denotes the normalised Lebesgue area measure on the open unit disc $\mathbb{D}$. 
The space  $D(\mu)$ is a reproducing kernel Hilbert space and   multiplication by the coordinate function $z$ turns out to be a bounded linear operator (refer \cite{aleman1993,richter1991}) on $D(\mu)$. For an elaborate discussion on Dirichlet spaces the reader may also consult \cite{dirichlet_sp2014}. S. Richter described a model for a cyclic, analytic $2$-isometry\cite{richter1991}. Indeed, a  cyclic, analytic $2$-isometry is unitarily equivalent to a multiplication operator $M_z$ on a Dirichlet space $D(\mu)$ for some finite, positive, Borel measure $\mu$ on unit circle $\mathbb{T}$.

D. Sarason \cite{sarason97} initiated the study of Dirichlet spaces by way of identifying such a space with a de Branges-Rovnyak space. This association has been further strengthened in recent years in the works of  \cite{cgr2022, cgr2010, kellay2015}. The identification of a Dirichlet space as a de Branges-Rovnyak space allows one to compute the reproducing kernel for the Dirichlet space. \\ Here, we include a brief description of de Branges-Rovnyak space for a ready reference. \\

\item {\bf de Branges-Rovnyak spaces:}
\\
For complex, separable Hilbert spaces $\mathcal{U}$, $\mathcal{V}$, let $\mathcal{B(U,V)}$ denote the Banach space of all bounded linear transformations from $\mathcal{U}$ to $\mathcal{V}$. The {\it Schur class} $S(\mathcal{U,V})$ is given by
\[S({\mathcal{U,V}})=\Big\{B:\mathbb{D} \to {\mathcal{B(U,V)}} : B \text{ is holomorphic, } \displaystyle\sup_{z\in \mathbb{D}}\|B(z)\|_{B({\mathcal{U,V}})}\leq 1\Big\}.\]
Observe that when ${\mathcal{U}}=\mathbb{C}=\mathcal{V}$ then the Schur class is nothing but the closed unit ball of $H^\infty(\mathbb{D})$, the set of all bounded holomorphic functions on $\mathbb{D}$. For any $B\in S(\mathcal{U,V})$, the de-Branges-Rovnyak space, $H(B)$ is the reproducing kernel Hilbert space associated with the $\mathcal{B(V)}$-valued semidefinite kernel given by
\[K_B(z,w)=\displaystyle\frac{I_{\mathcal{V}}-B(z)B(w)^*}{1-z\overline{w}},~~z,w\in \mathbb{D}.\]
For equivalent formulations of de Branges-Rovnyak spaces, the reader is referred to an excellent article by J. Ball \cite{Ball2015}. 
The kernel $K_B$ is normalised if $K_B(z,0)=I_{\mathcal{V}}$ for every $z\in \mathbb{D}$. This is equivalent to the condition $B(0)=0$. Further, when ${\mathcal{U}}=\mathbb{C}=\mathcal{V}$, we denote $H(B)$ by $H(b)$, the classical de Branges-Rovnyak space (Refer \cite{fricain_mashreghi_2016vol2} for the basic theory of the classical de Branges-Rovnyak spaces). 
\end{enumerate}

Let $\mathcal{H}$ be a complex, infinite dimensional, separable Hilbert space and $\mathcal{B}(\mathcal{H})$ denotes the $C^*$-algebra of bounded linear operators on $\mathcal{H}$. 
Here, we record a few definitions which are required in the sequel. 

Let $T\in {\mathcal{B(H)}}$.
We say that $T$ is {\it cyclic} if there exists a vector $f\in \mathcal{H}$ such that ${\mathcal{H}}=\bigvee \big\{T^nf : n\geq 0\big\}$ (vector $f$ is known as {\it cyclic vector}). An operator $T$ is said to be {\it analytic} if $\cap_{n\geq 0}T^n{\mathcal{H}}=\{0\}$. The  {\it Cauchy dual} $T'$ of a left invertible operator $T\in \mathcal{ B(H)}$ is defined as $T'=T(T^*T)^{-1}$. Following Agler \cite{agler1995_1}, an operator $T\in \mathcal{ B(H)}$ is said to be a {\it $2$-isometry} if \[I-2T^*T+T^{*2}T^2=0.\]
An operator $T\in \mathcal{B(H)}$ is said to be {\it subnormal} if there exist a Hilbert space $\mathcal{K}$ and an operator $S\in \mathcal{B(K)}$ such that $\mathcal{H} \subseteq K$, $S$ is normal and $S|_{\mathcal{H}}=T$. Readers are encouraged to refer \cite{conway1991} for a detailed study of subnormal operators. Agler proved in \cite{agler1985} that $T\in \mathcal{B(H)}$ is a subnormal contraction (that is, $T$ is subnormal and $\|T\|\leq 1$) if and only if 
\[B_n(T)\equiv \sum_{k=0}^{n}(-1)^k \binom{n}{k}T^{*k}T^k\geq 0 \text{ for all integers }n\geq 0.\]
Following \cite{At}, an operator $T\in B(\mathcal{H})$ is said to be {\it completely hyperexpansive} if \[B_n(T)\leq 0~~ \text{for all integers }n\geq 1.\]
Now we are ready to state and describe the problem that has been tackled in this article. 
\subsection{The Problem} 
We begin our discussion with the fact \cite[Theorem 6.1]{cgr2022} which states that for a cyclic, analytic $2$-isometry $T$ in $\mathcal{B(H)}$, the rank of $T^*T-I$ is finite if and only if there exist a finitely supported measure $\mu$ on the unit circle $\mathbb{T}$  such that $T$ is unitarily equivalent to a multiplication by the coordinate function $z$ on $D(\mu)$. In fact, if the rank of $T^*T-I$ is a positive integer $k,$ then the corresponding measure $\mu$ is supported at  exactly $k$  points of unit circle $\mathbb{T}.$
In view of this result, a construction of a cyclic, analytic 2-isometry reduces to choosing finitely many points on the unit circle and looking at the multiplication operator $M_z $ on the Dirichlet space $D(\mu)$, where $\mu$ is supported at the chosen points.

\noindent We now state the problem under consideration: 

\begin{problem}\label{p1}
Characterize  finitely supported, positive, Borel measures $\mu$ on unit circle $\mathbb{T}$ such that the Cauchy dual of $M_z$ on $D(\mu)$ is subnormal.
\end{problem}
\noindent An affirmative solution to the problem in the  case when $\mu$ is supported at a single point is given in \cite[Corollary 3.6]{badea2019} and \cite[Corollary 5.4]{cgr2022}.
Further, it is proved in \cite[Theorem 2.4]{cgr2022} that for a positive, Borel measure $\mu$ supported at any two antipodal points on unit circle $\mathbb{T}$, the Cauchy dual $M_z'$ of $M_z$ on $D(\mu)$ is subnormal. A counter-example to the CDSP has been constructed in \cite{mkvms2025} by choosing $\mu$ to be the sum of unit point mass measures at $\{1,i\}$, that is $\mu=\delta_1+\delta_i$. This result is further generalized in \cite{khasnis2025}, by choosing $\mu$ to be sum of unit point mass measures at $\{1,\zeta\}$ for arbitrary $\zeta\neq -1$ on the unit circle.\\
\noindent As a step towards the solution of Problem \ref{p1}, it is necessary to find a suitable generalization of the method of choosing $k$ points on the unit circle, for $k>2$. In view of the case when $k=2$, it is natural to choose equi-spaced points on the unit circle, and check the subnormality of the Cauchy dual $M_z'$ of $M_z$ on $D(\mu)$, where $\mu$ is supported at the chosen points. In the present paper, we deal with the case when $k=3$. In contrast with the case $k=2$, it turns out that the Cauchy dual $M_z'$ of $M_z$ on $D(\mu)$ is not subnormal in this case.
Indeed, the main result of this paper is as given below.
\begin{theorem}\label{main}
Let $\zeta_1$, $\zeta_2$ and $\zeta_3$ be three equi-spaced points on the unit circle $\mathbb{T}$. 
If the measure $\mu$ on the unit circle $\mathbb{T}$ is of the form  $\mu=\delta_{\zeta_1}+\delta_{\zeta_2}+\delta_{\zeta_3},$ then the Cauchy dual $M_z'$ of $M_z$ on the Dirichlet space $D(\mu)$ is not subnormal.
\end{theorem}

\noindent The rest of the paper is devoted to the proof of this result. 

\section{Proof of the main result}
Note that, in view of the rotational symmetry (see \cite[Proposition 7.1]{cgr2022}), it is sufficient to prove the main result by choosing $\zeta_1=1$, $\zeta_2=w$ and $\zeta_3=w^2,$ where $w$ denotes the complex cube-root of unity.

\noindent As mentioned earlier, the proof of the theorem relies on \cite[Theorem 6.4]{cgr2022}. The application of this theorem demands the discussion of the techniques used in the work of Costara \cite{costara2016}. In particular, the computations involved in obtaining the reproducing kernel of $D(\mu)$ in this case turn out to be crucial in what follows. For a ready reference, we briefly outline the work carried out in \cite{costara2016} and state the results that are used in the sequel.\\
Let $n\in \mathbb{N}$. Let $c_1,c_2,\dots,c_n$ be positive real numbers and $\zeta_1,\zeta_2,\dots,\zeta_n$ be distinct points on the unit circle $\mathbb{T}$. Consider the finite, positive, Borel measure  $\mu=\sum_{j=1}^{n}c_j\delta_{\zeta_j}$ on $\mathbb{T}.$
\begin{enumerate}
\item As an application of Riesz-Fejer theorem to the trigonometric polynomial\\ $\prod_{j=1}^{n}|z-\zeta_j|^2+\sum_{j=1}^{n}c_j\prod_{\substack{i=1\\i\neq j}}^{n}|z-\zeta_i|^2,$ there exist $\alpha_1,\alpha_2,\dots,\alpha_n \in \mathbb{C}\setminus\overline{\mathbb{D}}$ and $d>0$ such that the following equation holds:
\begin{equation}\label{eq3}
	\prod_{j=1}^{n}|z-\zeta_j|^2+\sum_{j=1}^{n}c_j\prod_{\substack{i=1\\i\neq j}}^{n}|z-\zeta_i|^2=d\prod_{j=1}^n |z-\alpha_j|^2,~~~z\in \mathbb{T}.
\end{equation}
\item The values of $\alpha_1,\alpha_2,\dots,\alpha_n$ as obtained in equation (\ref{eq3}) turn out to be important in the expressions for reproducing kernels of following two subspaces of $H^2$ obtained in \cite[Theorem 3.1, Theorem 4.4]{costara2016}:
	 Define \begin{equation}\label{eq16}
		O_\mu(z)=\displaystyle\frac{p(z)}{q(z)}
	\end{equation}
	with $p(z)=\displaystyle\frac{e^{i\theta}}{\sqrt{d}}\displaystyle\prod_{j=1}^{n}(z-\zeta_j)$ ($\theta \in \mathbb{R}$ is chosen such that $O_\mu(0)>0$) and\\ $q(z)=\displaystyle\prod_{j=1}^{n}(z-\alpha_j)$.\\
	\noindent We observe that $O_\mu$ is a rational function of degree $n$ having poles outside the closed unit disk, and simple zeros at $\zeta_j\in \mathbb{T}$, for $j=1,\dots,n$.  Then $\prod_{j=1}^{n}(z-\zeta_j)H^2=O_\mu H^2$ and hence $O_\mu H^2$ is a closed subspace of $D(\mu)$.
	The reproducing kernel for $O_\mu H^2$ is given by (Refer \cite[Theorem 3.1]{costara2016}) 
	$		\tilde{K_\mu}(z,\lambda)=\displaystyle \frac{O_\mu(z)}{1-\overline{\lambda}z}\overline{O_\mu(\lambda)}.
	$
	\item The reproducing kernel for $(O_\mu H^2)^\perp$ as a subspace of $D(\mu)$ is given by\\$
		\hat{K_\mu}(z,\lambda)=\displaystyle\sum_{j=1}^{n}g_j(\lambda)f_j(z) ~~~~ z,\lambda \in \mathbb{D}
	$
	where \begin{equation}\label{eq11}
		f_j(z)=\displaystyle\frac{O_\mu(z)}{O_\mu '(\zeta_j)(z-\zeta_j)}
	\end{equation}
	and 
	\begin{equation}\label{eq7}
		\begin{pmatrix}
			\overline{g_1(\lambda)}\\
			\overline{g_2(\lambda)}\\
			\vdots \\
			\overline{g_n(\lambda)}
		\end{pmatrix}=\left(\left[\langle f_i,f_j\rangle_\mu\right]_{1\leq i,j\leq n}\right)^{-1}
		\begin{pmatrix}
			f_1(\lambda)\\
			f_2(\lambda)\\
			\vdots \\
			f_n(\lambda)
		\end{pmatrix}.
	\end{equation}	
	\item The results \cite[Lemma 4.2, Lemma 4.3]{costara2016} provide formulae for obtaining $\langle f_i,f_j\rangle_\mu$ for any $1\leq i,j\leq n$ as follows :
	\begin{equation}\label{eq8}
		\|f_i\|_\mu^2=\langle f_i,f_i\rangle_\mu=c_i\zeta_i f_i'(\zeta_i)
	\end{equation} and for $i\neq j$,
	\begin{equation}\label{eq9}
		\langle f_i,f_j\rangle_\mu= \displaystyle\frac{1}{O_\mu'(\zeta_i)\overline{O_\mu'(\zeta_j)}(1-\zeta_i\overline{\zeta_j})}.
	\end{equation}
\end{enumerate}


\noindent As a first step towards the proof of the main theorem, we proceed with the following lemma. The lemma invokes the Riesz-Fejer theorem and derives a special case of equation (\ref{eq3}) by choosing $n=3,~ \zeta_1=1,~ \zeta_2=w, ~\zeta_3=w^2 \mbox{ and } c_1=c_2=c_3=1.$ 
\begin{lemma}\label{l1}
	For all $z\in \mathbb{T}$, there exists a polynomial $q(z)$ of degree $3$ and a constant $d>0$ such that \begin{equation}\label{eq13}
		|z-1|^2|z-w|^2|z-w^2|^2+|z-1|^2|z-w|^2+|z-1|^2|z-w^2|^2+|z-w|^2|z-w^2|^2=d~ |q(z)|^2.
	\end{equation}
\end{lemma}
\begin{proof}
	We have,
	\begin{align*}
		&|z-1|^2|z-w|^2|z-w^2|^2+|z-1|^2|z-w|^2+|z-1|^2|z-w^2|^2+|z-w|^2|z-w^2|^2\\
		&=(z-1)(\bar{z}-1)(z-w)(\bar{z}-\bar{w})(z-w^2)(\bar{z}-\bar{w^2})+(z-1)(\bar{z}-1)(z-w)(\bar{z}-\bar{w})\\
		&+(z-1)(\bar{z}-1)(z-w^2)(\bar{z}-\bar{w^2})+(z-w)(\bar{z}-\bar{w})(z-w^2)(\bar{z}-\bar{w^2})\\
		&=11-6(1+\bar{w}+\bar{w^2})z-6(1+w+w^2)\bar{z}+3(1+\bar{w}+\bar{w^2})z^2+3(1+w+w^2)\bar{z}^2-z^3-\bar{z}^3\\
		&=11-z^3-\left(\frac{1}{z}\right)^3\\
		&=-\displaystyle\frac{1}{z^3}\left[z^6-11z^3+1\right]
	\end{align*}
	\noindent Let $f(z)=z^6-11z^3+1$. Note that the polynomial $f$ has a property that for $z_0\in \mathbb{C}$ if $f(z_0)=0$, then $f\left(\displaystyle\frac{1}{\overline{z_0}}\right)=0$. Therefore three of the roots of $f(z)$ lie outside the closed unit disc. We denote these roots by $\alpha_1, \alpha_2$ and $\alpha_3$ and the other three roots by $\beta_1$, $\beta_2$ and $\beta_3$ where $\beta_1=\displaystyle\frac{1}{\overline{\alpha_1}}$,  $\beta_2=\displaystyle\frac{1}{\overline{\alpha_2}}$ and $\beta_3=\displaystyle\frac{1}{\overline{\alpha_3}}$. Thus, the roots of $f(z)$ are $\alpha_1, \alpha_2, \alpha_3 \displaystyle\frac{1}{\overline{\alpha_1}}, \displaystyle\frac{1}{\overline{\alpha_2}}$ and $\displaystyle\frac{1}{\overline{\alpha_3}}$.\\
	\noindent At this point, we observe that $f(z)$ is a quadratic in $z^3$, which implies that\\ $z^3=\displaystyle\frac{11\pm3\sqrt{13}}{2}$. Let $\alpha=\left(\displaystyle\frac{11+3\sqrt{13}}{2}\right)^{\frac{1}{3}}\in\mathbb{C}\backslash\overline{\mathbb{D}}$,
	we obtain the expression $z^3-\alpha^3=0$. This implies that $\alpha_1=\alpha, \alpha_2=\alpha w$ and $\alpha_3=\alpha w^2$.\\
	
	

	\noindent We now define $q(z)$ as $q(z)=(z-\alpha_1)(z-\alpha_2)(z-\alpha_3)=z^3-\alpha^3$. \\
	Let $\beta=\frac{1}{\alpha}$ then, for $z\in \mathbb{T}$, we consider
	\begin{align*}
		\displaystyle\frac{(z-\beta_1)(z-\beta_2)(z-\beta_3)}{z^3}&=1-\displaystyle\frac{\beta^3}{z^3}=1-\displaystyle\frac{1}{\alpha^3 z^3}\\
		&=\displaystyle\frac{1}{\alpha^3}(\alpha^3-\bar{z}^3)\\
		&=\displaystyle-\frac{1}{\alpha^3}(\overline{q(z)})
	\end{align*}
 Note that since left side of equation (\ref{eq13}) is positive, we have, $-\displaystyle\frac{1}{z^3}f(z)>0$, for all $z\in \mathbb{T}$.\\
	Let $b=\alpha^3$ then we get, \[-\displaystyle\frac{1}{z^3}f(z)=\frac{1}{b}q(z)\overline{q(z)}=\frac{1}{b}\big|q(z)\big|^2,\text{ for all }z\in \mathbb{T}\]	
	This proves that $b>0$. Thus we get, 
	\[|z-1|^2|z-w|^2|z-w^2|^2+|z-1|^2|z-w|^2+|z-1|^2|z-w^2|^2+|z-w|^2|z-w^2|^2=d~ |q(z)|^2\] where $d=\displaystyle\frac{1}{b}>0$.
\end{proof}

\setcounter{equation}{11}
\noindent The computation of $q(z)$ allows us to apply formula in equation (\ref{eq16}), as applicable in the present case. We thus have the following expression for $O_\mu(z)$: \begin{equation}\label{eq1}
	O_\mu(z)=\displaystyle\frac{(z^3-1)}{\sqrt{d}q(z)}.
\end{equation}
Observe that we take $\theta=0$ so that $O_\mu(0)>0$.

\noindent Now,
$$O_\mu'(z)=\displaystyle\frac{3z^2(1-\alpha^3)}{\sqrt{d}(z^3-\alpha^3)^2}$$ which implies that $$O_\mu'(1)=\displaystyle\frac{3}{\sqrt{d}(1-\alpha^3)},~ O_\mu'(w)=\displaystyle\frac{3w^2}{\sqrt{d}(1-\alpha^3)}=w^2O_\mu'(1),~ O_\mu'(w^2)=\displaystyle\frac{3w}{\sqrt{d}(1-\alpha^3)}=wO_\mu'(1)$$ Using the fact that $\alpha$ is a root of $f(z)$, we obtain $|O_\mu'(1)|=|O_\mu'(w)|=|O_\mu'(w^2)|=1$.\\
\noindent	We now use equation (\ref{eq11}) to obtain the following expressions for $f_1(z)$, $f_2(z)$ and $f_3(z)$:\\
	
	$f_1(z)=\displaystyle\frac{O_\mu(z)}{O_\mu'(1)(z-1)}=\displaystyle\frac{(1-\alpha^3)(z-w)(z-w^2)}{3(z^3-\alpha^3)}~$, $~f_2(z)=\displaystyle\frac{O_\mu(z)}{O_\mu'(w)(z-w)}=\displaystyle\frac{(1-\alpha^3)(z-1)(z-w^2)}{3(z^3-\alpha^3)}$   and 
	$~f_3(z)=\displaystyle\frac{O_\mu(z)}{O_\mu'(w^2)(z-w^2)}=\displaystyle\frac{(1-\alpha^3)(z-1)(z-w)}{3(z^3-\alpha^3)}.$ \\
	
	\noindent

	\noindent Now, equations (\ref{eq8}) and (\ref{eq9}) are used, to obtain the matrix: 
	$B=D^{-1}$ where $D=((\langle f_i,f_j\rangle ))_{1\leq i,j\leq 3}.$\\
	A simple computation reveals that,\\
	$\|f_1\|^2=\|f_2\|^2=\|f_3\|^2=-\displaystyle\frac{(2+\alpha^3)}{(1-\alpha^3)}=-\displaystyle\frac{(2+b)}{(1-b)}$.\\
	Using the fact that $\alpha$ is a root of $f(z),$ we have \\
	$\langle f_1,f_2\rangle=\langle f_2,f_3\rangle=\displaystyle\frac{(1-\alpha^3)^2}{9\alpha^3(w-1)}=\displaystyle\frac{1}{(w-1)}$ and $\langle f_1,f_3\rangle=\displaystyle\frac{(1-\alpha^3)^2}{9\alpha^3(w^2-1)}=\displaystyle\frac{1}{(w^2-1)}.$\\
	Hence 
	$D=\begin{bmatrix}
		-\displaystyle\frac{(2+b)}{(1-b)} & \displaystyle\frac{1}{(w-1)} & \displaystyle\frac{1}{(w^2-1)}\\
		\displaystyle\frac{1}{(w^2-1)} & -\displaystyle\frac{(2+b)}{(1-b)} & \displaystyle\frac{1}{(w-1)}\\
		\displaystyle\frac{1}{(w-1)} & \displaystyle\frac{1}{(w^2-1)} & -\displaystyle\frac{(2+b)}{(1-b)} \\
	\end{bmatrix}.$\\
This matrix can be written in the form
$D=\begin{bmatrix}
	x & s & \bar{s}\\
\bar{s} & x & s\\
	s & \bar{s} & x\\
\end{bmatrix},$ where $x=-\displaystyle\frac{(2+b)}{(1-b)}$ and $s=\displaystyle\frac{1}{(w-1)}$.
Further, we have $\det D=x(x^2-1)$.
Therefore, we get
	\begin{equation}\label{mat1}
		B=D^{-1}=\displaystyle\frac{1}{x(x^2-1)}\begin{bmatrix}
	{(x^2-\frac{1}{3})} & {\bar{s}^2-xs} & {s^2-x\bar{s}}\\
{s^2-x\bar{s}} & {(x^2-\frac{1}{3})} & {\bar{s}^2-xs} \\
	{\bar{s}^2-xs} & {s^2-x\bar{s}} & {(x^2-\frac{1}{3})} \\
\end{bmatrix}.\end{equation}

\noindent Before we present the proof of the main theorem, we quote the following results from \cite{cgr2022} which are used in the proof of the theorem.
\begin{theorem}\cite[Theorem 6.4]{cgr2022}\label{thm2}
	For positive scalars $c_1,c_2,\dots,c_k$ and distinct points $\zeta_1,\zeta_2,\dots,\zeta_k$ on the unit circle $\mathbb{T}$, consider the positive Borel measure $\mu=\sum_{j=1}^{k}c_j\delta_{\zeta_j}$ on $\mathbb{T}$, where $\delta_{\zeta_j}$ denotes the Dirac delta measure supported at $\zeta_j$. Let $X(z)=(z,z^2,\dots,z^k)^T$ and $\{e_j\}_{j=1}^k$ denote the standard basis of $\mathbb{C}^k$. Then there exist $\alpha_1, \alpha_2,\dots,\alpha_k\in \mathbb{C}\setminus\overline{\mathbb{D}}$ and a $k\times k$ upper triangular matrix $P$ such that the Dirichlet space $D(\mu)$ coincides with the de Branges-Rovnyak space $H(B)$ with equality of norms, where $B=\left(\frac{p_1}{q}, \frac{p_2}{q}, \dots, \frac{p_k}{q}\right)$ and 
	\[p_j(z)=\left\langle PX(z),e_j\right\rangle~,~j=1,\dots,k,~~~q(z)=\prod_{j=1}^k (z-\alpha_j)\]
	Moreover, $\alpha_1,\dots,\alpha_k$ are governed by
	\[\prod_{j=1}^{k}|z-\zeta_j|^2+\sum_{j=1}^{k}c_j\prod_{\substack{i=1\\i\neq j}}^{k}|z-\zeta_i|^2=d\prod_{j=1}^k |z-\alpha_j|^2,~~~z\in \mathbb{T}\]
	for some $d>0$.\\
\end{theorem}

\begin{theorem}\cite[Theorem 2.1]{cgr2022}\label{thm3}
	Let $B=(b_1,\dots,b_k)\in {\mathcal{S}}(\mathbb{C}^k,\mathbb{C})$ be such that $B(0)=0$ where
	\[b_j(z)=\displaystyle\frac{p_j(z)}{\prod_{j=1}^k (z-\alpha_j)}\]
	for polynomials $p_j$ of degree at most $k$ and distinct numbers $\alpha_1,\dots,\alpha_k \in \mathbb{C}\setminus\overline{\mathbb{D}}$. for $r=1,\dots, k$, let $a_r=\prod_{1\leq t\neq r \leq k} (\alpha_r-\alpha_t)$. Assume that the operator $M_z$ of multiplication by $z$ on the de Branges-Rovnyak space $H(B)$ is bounded. Then the Cauchy dual $M_z'$ of $M_z$ is subnormal if and only if the matrix \\
	\begin{center}{$\displaystyle\sum_{r,t=1}^{k}\left( \frac{1}{a_r\overline{a_t}}\sum_{j=1}^{k}p_j(\alpha_r)\overline{p_j(\alpha_t)}\right)\left(1-\frac{1}{\alpha_r \overline{\alpha_t}}\right)^l\left(\left(\frac{1}{\alpha_r^{m+2}\overline{\alpha_t}^{n+2}}\right)\right)_{m,n\geq 0}$}
		\end{center}
	is formally positive semi-definite for every $l\geq 1$.
\end{theorem}
\begin{cor}\cite[Corollary 4.3]{cgr2022}\label{cor1}
	Assume the hypothesis of Theorem \ref{thm3}. If $\alpha_r\overline{\alpha_t}\notin [1,\infty)$ for every $1\leq r\neq t\leq k$, then $M_z'$ is subnormal if and only if \\
	$\sum_{j=1}^k p_j(\alpha_r)\overline{p_j(\alpha_t)}=0, ~1\leq r\neq t\leq k.$
\end{cor}
\noindent We now proceed to give the proof of the main theorem.
\begin{proof}
	Note that by Theorem \ref{thm2}, we conclude that $D(\mu)=\mathcal{H}(B)$ with equality of norms where $B=\left(\displaystyle \frac{p_1}{q},\frac{p_2}{q},\frac{p_3}{q}\right).$ In order to reach the desired conclusion, appealing to Corollary \ref{cor1}, we now show that $\displaystyle \sum_{j=1}^3 p_j(\alpha_r)\overline{p_j(\alpha_t)}\neq 0, ~1\leq r\neq t\leq 3.$\\
\noindent By \cite[Remark 6.5]{cgr2022} we have,
\begin{align*}
		\sum_{j=1}^{3}p_j(z)\overline{p_j(u)} &=q(z)\overline{q(u)}-p(z)\overline{p(u)}\left(1+\sum_{i,j=1}^{3}\displaystyle\frac{\overline{b_{ji}}}{O_\mu'(\zeta_j)\overline{O_\mu'(\zeta_i)}}\frac{(1-z\bar{u})}{(z-\zeta_j)(\bar{u}-\bar{\zeta_i})}\right)\\
		&=q(z)\overline{q(u)}-p(z)\overline{p(u)}-p(z)\overline{p(u)}(1-z\bar{u})\left(\sum_{i,j=1}^{3}\displaystyle\frac{\overline{b_{ji}}}{O_\mu'(\zeta_j)\overline{O_\mu'(\zeta_i)}}\frac{1}{(z-\zeta_j)(\bar{u}-\bar{\zeta_i})}\right).
	\end{align*}
	
		The numbers $b_{ij}$, the entries of the matrix $B$, can be read from the expression (\ref{mat1}). Thus we have $b_{11}=b_{22}=b_{33}=\displaystyle\frac{1}{x(x^2-1)}{(x^2-\frac{1}{3})}$,\\ $b_{12}=b_{23}=b_{31}=\displaystyle\frac{1}{x(x^2-1)}{(\bar{s}^2-xs)} $ and
$b_{13}=b_{21}=b_{32}=\displaystyle\frac{1}{x(x^2-1)}{(s^2-x\bar{s})}.$\\
	
	Consider the expression,
	\begin{align*}
		\sum_{i,j=1}^{3}\displaystyle\frac{\overline{b_{ji}}}{O_\mu'(\zeta_j)\overline{O_\mu'(\zeta_i)}}\frac{1}{(z-\zeta_j)(\bar{u}-\bar{\zeta_i})}=& \displaystyle{\frac{\overline{b_{11}}}{|O_\mu'(1)|^2}\frac{1}{(z-1)(\bar{u}-1)}+\frac{\overline{b_{12}}}{O_\mu'(1)\overline{O_\mu'(w)}}\frac{1}{(z-1)(\bar{u}-\overline{w})}}\\
		&\displaystyle{+\frac{\overline{b_{13}}}{O_\mu'(1)\overline{O_\mu'(w^2)}}\frac{1}{(z-1)(\bar{u}-\overline{w^2})}		+\frac{\overline{b_{21}}}{O_\mu'(w)\overline{O_\mu'(1)}}\frac{1}{(z-w)(\bar{u}-{1})}}\\
		&+\displaystyle{\frac{\overline{b_{22}}}{|O_\mu'(w)|^2}\frac{1}{(z-w)(\bar{u}-\overline{w})}+\frac{\overline{b_{23}}}{O_\mu'(w)\overline{O_\mu'(w^2)}}\frac{1}{(z-w)(\bar{u}-\overline{w^2})}}\\
		&+\displaystyle{\frac{\overline{b_{31}}}{O_\mu'(w^2)\overline{O_\mu'(1)}}\frac{1}{(z-w^2)(\bar{u}-{1})}+\frac{\overline{b_{32}}}{O_\mu'(w^2)\overline{O_\mu'(w)}}\frac{1}{(z-w^2)(\bar{u}-\overline{w})}}\\
		&+\displaystyle{\frac{\overline{b_{33}}}{|O_\mu'(w^2)|^2}\frac{1}{(z-w^2)(\bar{u}-\overline{w^2})}}\\
		=&\overline{b_{11}}\left[\displaystyle\frac{1}{(z-1)(\bar{u}-1)}+\frac{1}{(z-w)(\bar{u}-\overline{w})}+\frac{1}{(z-w^2)(\bar{u}-\overline{w^2})}\right]\\
		&+\displaystyle\frac{\overline{b_{12}}}{w}\left[\frac{1}{(z-1)(\bar{u}-\bar{w})}+\frac{1}{(z-w)(\bar{u}-\overline{w^2})}+\frac{1}{(z-w^2)(\bar{u}-1)}\right]\\
		&+\displaystyle\frac{{b_{12}}}{w^2}\left[\displaystyle\frac{1}{(z-1)(\bar{u}-\overline{w^2})}+\frac{1}{(z-w)(\bar{u}-1)}+\frac{1}{(z-w^2)(\bar{u}-\overline{w})}\right]\\
	=	&\displaystyle\frac{1}{(z^3-1)(\bar{u}^3-1)}\left[3\left(\overline{b_{11}}+\frac{\overline{b_{12}}}{w}+\frac{{b_{12}}}{w^2}\right)(z\bar{u})^2\right.\\
		&\left.\displaystyle+3\left(\overline{b_{11}}+\frac{\overline{b_{12}}}{w^2}+\frac{{b_{12}}}{w}\right)(z\bar{u})+3(\overline{b_{11}}+\overline{b_{12}}+b_{12})\right].
	\end{align*}
	Simplifying, we get,
	\[\sum_{j=1}^{3}p_j(z)\overline{p_j(u)}=\displaystyle\left((1-b)+\frac{3b}{x+1}\right)(z\bar{u})^3+\frac{3b}{x(x+1)}(z\bar{u})^2+\frac{3b}{x(x-1)}(z\bar{u}).\]
	We now choose $z=\alpha, u=\alpha w$ (two of the three roots of $f(z)$ lying outside the closed unit disk), and substitute these values in the above equation. Thus
	\[\sum_{j=1}^{3}p_j(\alpha)\overline{p_j(\alpha w)}=\displaystyle(\alpha^2\overline{w})\left[\left((1-b)+\frac{3b}{x+1}\right)(\alpha^2\overline{w})^2+\frac{3b}{x(x+1)}(\alpha^2\overline{w})+\frac{3b}{x(x-1)}\right].\]
	It can be seen that $(\alpha^2\overline{w})$ is not a root of the quadratic equation\\ $\left((1-b)+\frac{3b}{x+1}\right)Y^2+\frac{3b}{x(x+1)}Y+\frac{3b}{x(x-1)}=0$ implying that
	$\sum_{j=1}^{3}p_j(\alpha)\overline{p_j(\alpha w)}\neq0.$
	
	Therefore, the Cauchy dual $M_z'$ of $M_z$ on $D(\mu)$ is not subnormal.
\end{proof}



\section{Epilogue}
The authors believe that the main theorem of the present article can be proved for the measure $\mu$ in the general form, i.e. $\mu= c_1\delta_{\zeta_1}+c_2\delta_{\zeta_2}+c_3\delta_{\zeta_3},$ where $c_1,c_2,c_3$ are positive real numbers, but the computational complexity increases substantially.  Further, there is a strong evidence that the following conjecture has an affirmative solution.
\begin{conj}
If $\zeta_1,\zeta_2,\zeta_3$ are any three distinct points on the unit circle and if $\mu=c_1\delta_{\zeta_1}+c_2\delta_{\zeta_2}+c_3\delta_{\zeta_3}$, then the Cauchy dual $M_z'$ of $M_z$ on $D(\mu)$ is not subnormal.
\end{conj}
\section*{Declarations}

The present work is carried out at the research center at the Department of Mathematics, S. P. College, Pune, India(autonomous).






\bibliographystyle{plain}
\bibliography{ref.bib}


\end{document}